\documentclass{article}
\usepackage{amsfonts}
\usepackage{amsmath}

\newtheorem{theorem}{Theorem}

\newtheorem{definition}[theorem]{Definition}

\newtheorem{proposition}[theorem]{Proposition}

\input{tcilatex}
\begin{document}

\title{The homogeneous flow of a \\
parallelizable manifold}
\author{Erc\"{u}ment Orta\c{c}gil}
\maketitle

\begin{abstract}
Motivated by the Hamilton's Ricci flow, we define the homogeneous flow of a
parallelizable manifold and show the short time existence and uniqueness of
its solutions. We indicate the relation of this flow to the Poincare
Conjecture.
\end{abstract}

\QTP{Body Math}
\bigskip $\bigskip $\textbf{Contents}

\QTP{Body Math}
1) Introduction

\QTP{Body Math}
2) Parallelizable manifolds and local Lie groups

3) The homogeneous flow

4) Gauge group

5) Poincare Conjecture

6) Appendix: The relation of HF to the Ricci flow

\section{Introduction}

In [8] we proposed a generalization of the Klein's Erlangen Program based on
our earlier work in arXiv. A parallel theory is proposed in [2]. This note
is the result of our efforts to give an alternative approach to the proof of
the Poincare Conjecture using a new geometric evolution equation which
emerges out of this program.

The central concept in the framework of [8] is that of a prehomogeneous
geometry (phg) and its curvature. The order of a phg is the order of jets
involved in its definition. The curvature is the obstruction to the local
homogeneity of the phg. In this note we are interested in the simplest phg
of order zero, i.e, a parallelizable manifold $(M,\varepsilon )$ where $%
\varepsilon $ denotes the parallelization. If the curvature $\mathfrak{R}%
(\varepsilon )$ vanishes, $M$ becomes locally homogeneous in two ways and is
called a local Lie group in [1]. If $M$ is also simply connected and $%
\varepsilon $ is complete, then $M$ is the homogeneous space of two global
and simply transitive transformation groups which correspond to the
left-right actions of a Lie group. Section 2 contains a concise exposition
of this theory with more details than in [1] on certain points, also
clarifying certain ambiguities in [1]. It is worth stressing here that the
theory of local Lie groups is not a simple consequence of the present global
theory but has its own set of interesting and delicate geometric structures
as stated in [5] which deeply inspired our work. For instance, a local Lie
group in this sense does not always imbed in a global Lie group ([5]). In
fact, it is shown in [1] that the opposite is true: a Lie group is a special
(globalizable) local Lie group! Therefore, in the words of [6], Section 2
"reinstates the paradigm of local to global to its historical record".

In Section 3 we define the homogeneous flow (HF) of a parallelizable
manifold which is inspired by the Ricci flow of Hamilton. HF is a second
order nonlinear evolution equation which starts with an arbitrary
parallelism and flows towards a parallelism with vanishing curvature. We
show that HF is weakly parabolic. Using the DeTurck trick [3], we show that
HF is equivalent to a strongly parabolic flow thus establishing the
existence and uniqueness of the short time solutions of HF.

Using the first order universal gauge group which is an infinite dimensional
Frechet Lie group, we show in Section 3 that the short time solutions of HF
localize to solutions of an ODE defined at each point of $M.$ Therefore the
evolution of the initial parallelism at some point is determined for all
times by the value of the curvature at $t=0$ at that point, indicating the
simple nature of HF. However, the explicit form of this ODE also indicates
that HF will develop finite time singularities like the Ricci flow. The
nature of these singularities remains to be studied.

In Section 4 we briefly comment on the relation of HF to PC which has been
the main motivation for this note.

In the Appendix we explain briefly how we hit upon HF while trying to
understand the Ricci flow.

Finally, it is worth stressing here that HF is defined for any phg (in
particular for a Riemannian geometry as a phg of order one), the key fact
being that the top principal bundle defined by the phg is parallelizable
([8]). The curvature of a Riemannian geometry as a phg vanishes if and only
if the underlying metric has constant sectional curvature (which is
equivalent to local homogeneity. See page 6 of [2] for a simple formula for
this curvature).

\section{Parallelizable manifolds and local Lie groups}

Let $M$ be a smooth manifold with $\dim M\geq 2$ and $\mathcal{U}_{k}(M)$
(shortly $\mathcal{U}_{k})$ be the universal groupoid of order $k$ on $M.$
The elements of $\mathcal{U}_{k}$ are the $k$-jets of local diffeomorphisms
of $M.$ We call an element of $\mathcal{U}_{k}$ with source at $p$ and
target at $q$ a $k$-arrow from $p$ to $q$ and denote it by $j_{k}(f)^{p,q}.$
Therefore $\mathcal{U}_{0}=M\times M$ is the pair groupoid. The relevant
universal groupoids in this section are $\mathcal{U}_{0}$ and $\mathcal{U}%
_{1}$. The projection homomorphism $\pi :\mathcal{U}_{1}\rightarrow \mathcal{%
U}_{0}$ of groupoids maps a $1$-arrow from $p$ to $q$ to the pair $(p,q).$ A
splitting $\varepsilon :\mathcal{U}_{0}\rightarrow \mathcal{U}_{1}$ is a
homomorphism of groupoids so that $\pi \circ \varepsilon =id_{\mathcal{U}%
_{0}}.$ Thus $\varepsilon $ assigns to any pair $(p,q)$ a unique $1$-arrow
from $p$ to $q$ and this assignment preserves the composition and inversions
of arrows. We easily check that $\pi :\mathcal{U}_{1}\rightarrow \mathcal{U}%
_{0}$ admits a splitting if and only if $M$ is parallelizable. If $p\in
(U,x^{i})$ has coordinates $\overline{x}^{i}$ and $q\in (V,y^{i})$ has
coordinates $\overline{y}^{i}$, then $\varepsilon (p,q)$ has the local
representation $\varepsilon _{j}^{i}(\overline{x}^{1},...,\overline{x}^{n},%
\overline{y}^{1},...,\overline{y}^{n})=\varepsilon _{j}^{i}(\overline{x},%
\overline{y}),$ $1\leq i,j\leq n=\dim M.$ Thus we have the coordinate
formulas

\begin{eqnarray}
\varepsilon _{a}^{i}(z,y)\varepsilon _{j}^{a}(x,z) &=&\varepsilon
_{j}^{i}(x,y)  \notag \\
\varepsilon _{j}^{i}(x,x) &=&\delta _{j}^{i}  \notag \\
\varepsilon _{a}^{i}(y,x)\varepsilon _{j}^{a}(x,y) &=&\delta _{j}^{i}
\end{eqnarray}

We use summation convention in (1). In this section we fix the splitting $%
\varepsilon $ once and for all and let $(M,\varepsilon )$ denote the
parallelizable manifold $M.$

Now we consider the first order nonlinear PDE

\begin{equation}
\frac{\partial f^{i}(x)}{\partial x^{j}}=\varepsilon _{j}^{i}(x,f(x))
\end{equation}%
for some local diffeomorphism $y^{i}=f^{i}(x).$ The integrability conditions
of (2) are given by

\begin{equation}
\mathcal{R}_{jk}^{i}(x,y)\overset{def}{=}\left[ \frac{\partial \varepsilon
_{k}^{i}(x,y)}{\partial x^{j}}+\frac{\partial \varepsilon _{k}^{i}(x,y)}{%
\partial y^{a}}\varepsilon _{j}^{a}(x,y)\right] _{[jk]}=0
\end{equation}%
where $[jk]$ denotes the alternation of the indices $j,k.$ We have $\mathcal{%
R}(p,q)\in \wedge ^{2}T_{p}^{\ast }\otimes T_{q}.$ If (3) admits a solution $%
f$ with $f(p)=q$ for any $(p,q)\in U\times V$, then clearly $\mathcal{R}=0$
on $U\times V.$ Conversely, by the well known existence and uniqueness
theorem for the first order systems of PDE's, if $\mathcal{R}=0$ on $U\times
V$, then we may assign any pair $(p,q)\in U\times V$ as initial condition
and solve (2) uniquely for some $f$ defined on $\overline{U}\subset U$
satisfying $f(p)=q.$ Further, $j_{1}(f)^{x,f(x)}\in \varepsilon (\mathcal{U}%
_{0})$ for all $x\in \overline{U}$ and we may choose $\overline{U}=U$ if $U$
is simply connected. Note that $\mathcal{R(}p,p)=0$ for all $p\in M.$

\begin{definition}
$\mathcal{R}$ is the groupoid curvature of $(M,\varepsilon )$ and $%
(M,\varepsilon )$ is locally homogeneous (or a local Lie group) if $\mathcal{%
R}=0$ on $M\times M.$
\end{definition}

To justify the term local homogeneity, we assume $\mathcal{R}=0$ and let $%
\mathcal{S}$ denote the set of all local solutions of (2). Since $%
\varepsilon $ is a homomorphism of groupoids, $\mathcal{S}$ is easily seen
to be a pseudogroup. Some $f\in \mathcal{S}$ is determined on its domain by
any of its $0$-arrows $(p,f(p)).$ Now let $f\in \mathcal{S}$ be defined on $%
U $, $p\in U$ and $C$ a path from $p$ to some $q\in M.$ We can "analytically
continue" $f$ along $C$ but may not be able to "reach" $q.$ We call $%
(M,\varepsilon )$ complete if all elements of $\mathcal{S}$ can be continued
indefinitely along all paths in $M.$ Note that we define the completeness of 
$(M,\varepsilon )$ only when $\mathcal{R}=0$ (at least here). Assuming
completeness, two paths from $p$ to $q$ may give different values at $q$ if
these paths are not homotopic. However, the standard monodromy argument
shows that we get the same values at $q$ if these paths are homotopic. In
particular, if $M$ is simply connected, we easily see that any $f\in 
\mathcal{S}$ extends to a global diffeomorphism of $M.$ Further, these
global transformations are closed under composition and inversion and
therefore they form a global transformation group of $M$ which acts simply
transitively. We continue to denote this transformation group by $\mathcal{S}
$ and call $\mathcal{S}$ globalizable (as a pseudogeoup). Note that $%
\mathcal{S}$ may be globalizable without $M$ being simply connected but $%
\mathcal{S}$ is of course complete if it is globalizable and what we have
shown above is that completeness together with simple connectedness implies
globalizability. If $(M,\varepsilon )$ is complete but not globalizable,
then we can lift $\mathcal{S}$ to a pseudogroup $\mathcal{S}^{u}$ on the
universal cover $M^{u}$ of $M$ and globalize $\mathcal{S}^{u}$ on $M^{u}$
such that the covering transformations form a discontinuous subgroup of $%
\mathcal{S}^{u}$ isomorphic to the fundamental group of $M.$

If $\mathcal{R}=0$, there is another pseudogroup on $M$ defined as follows.
Let $f(a,b,z)$ denote the unique local solution of (2) in the variable $z$
satisfying the initial condition $a\rightarrow b.$ We fix some $p,q\in
(U,x^{i})$ and define

\begin{equation}
\widetilde{\varepsilon }_{j}^{i}(p,q)\overset{def}{=}\left( \frac{\partial
f^{i}(p,x,q)}{\partial x^{j}}\right) _{x=p}
\end{equation}

Note that $\widetilde{\varepsilon }(p,q)$ is defined for close $p,q$ unless $%
\mathcal{S}$ is globalizable. We check that $\widetilde{\varepsilon }$ is a 
\textit{local }splitting of $\pi :\mathcal{U}_{1}\rightarrow \mathcal{U}%
_{0}. $ Therefore we can replace (2) by

\begin{equation}
\frac{\partial h^{i}(x)}{\partial x^{j}}=\widetilde{\varepsilon }%
_{j}^{i}(x,h(x))
\end{equation}

Now the local diffeomorphism $h:x\rightarrow f(p,x,q)$ satisfies $h(p)=q$
and solves (5). In particular the integrability conditions of (5) are
satisfied. Thus we get a pseudogroup $\widetilde{\mathcal{S}}$ in the same
way we get $\mathcal{S}$. The only difference is that $\widetilde{\mathcal{S}%
}$ is locally transitive whereas $\mathcal{S}$ is globally transitive. The
elements of $\mathcal{S}$ and $\widetilde{\mathcal{S}}$ commute whenever
their compositions are defined. If $\mathcal{S}$ globalizes, then so does $%
\widetilde{S}$ in which case we get two global commuting transformation
groups of $M.$ Now it is easy to construct an abstract Lie group $G$ whose
underlying manifold is $M$ and its left-right (or right-left) translations
can be identified with $\mathcal{S}$ and $\widetilde{\mathcal{S}}.$ However,
note that there is no such canonical identification!

Up to now we assumed $\mathcal{R}=0$ and dealt with the "Lie group". Now we
drop the assumption $\mathcal{R}=0$ and consider the parallelizable manifold 
$(M,\varepsilon )$ our purpose being to construct the "Lie algebra".

We define

\begin{equation}
\Gamma _{jk}^{i}(x)\overset{def}{=}\left( \frac{\partial \varepsilon
_{k}^{i}(x,y)}{\partial y^{j}}\right) _{y=x}
\end{equation}

It is extremely crucial that $\Gamma _{jk}^{i}(x)$ need not be symmetric in $%
j,k.$ Differentiating the third formula in (1) with respect to $x$ at $y=x$
gives

\begin{equation}
\left( \frac{\partial \varepsilon _{k}^{i}(x,y)}{\partial x^{j}}\right)
_{y=x}=-\Gamma _{jk}^{i}(x)
\end{equation}

The $1$-arrow $\varepsilon (p,q)$ induces an isomorphism $\varepsilon
(p,q)_{\ast }$ of the tangent spaces $\varepsilon (p,q)_{\ast
}:T_{p}\rightarrow T_{q}$ which extends to an isomorphism $\varepsilon
(p,q)_{\ast }:\left( T_{r}^{m}\right) _{p}\rightarrow \left(
T_{r}^{m}\right) _{q}$ of the tensor spaces. A tensor field $t$ is $%
\varepsilon $-parallel if $\varepsilon (p,q)_{\ast }t(p)=t(q)$ for all $%
p,q\in M.$ Thus an $\varepsilon $-parallel $t$ is globally determined by its
value at any point. For instance, the tensor $(t_{j}^{i})$ is $\varepsilon $%
-parallel if and only if

\begin{equation}
t_{j}^{i}(x)=\varepsilon _{a}^{i}(p,x)t_{b}^{a}(p)\varepsilon _{j}^{b}(x,p)
\end{equation}%
for any fixed but arbitrary $p$ and all $x.$ Differentiating (8) with
respect to $x$ at $x=p,$ substituting from (6), (7) and omitting $p$ from
our notation, we get

\begin{equation}
\nabla _{r}t_{j}^{i}\overset{def}{=}\frac{\partial t_{j}^{i}}{\partial x^{r}}%
-\Gamma _{ra}^{i}t_{j}^{a}+\Gamma _{rj}^{a}t_{a}^{i}=0
\end{equation}

The operator $\nabla $ extends to all tensor fields in the obvious way. Note
that our sign convention in (9) is the opposite of the one used in tensor
calculus because of our choice of (6) rather than (7) but this point is not
important. It is crucial that $r$ is the first index in $\Gamma _{r\bullet
}^{\bullet }$ in (9). The derivation of (9) from (8) proves that $%
\varepsilon $-parallelity of $t$ implies $\nabla t=0.$ Converse is also
true. To see this, let $\xi =(\xi ^{i})$ be a vector field satisfying

\begin{equation}
\nabla _{r}\xi ^{i}=\frac{\partial \xi ^{i}}{\partial x^{r}}-\Gamma
_{ra}^{i}\xi ^{a}=0
\end{equation}

The integrability conditions of (10) are given by

\begin{equation}
\widetilde{\mathfrak{R}}_{rj,k}^{i}\overset{def}{=}\left[ \frac{\partial
\Gamma _{jk}^{i}}{\partial x^{r}}+\Gamma _{rk}^{a}\Gamma _{ja}^{i}\right]
_{[rj]}=0
\end{equation}

The order of the indices is quite relevant in (11). If $\widetilde{\mathfrak{%
R}}=0$ is identically satisfied on $M,$ then for any initial condition $\xi
^{i}(p)$ at some $p\in M,$ we have a unique solution $\xi ^{i}(x)$ of (11)
around $p$ satisfying this initial condition. However, $\xi ^{i}(p)$
determines an $\varepsilon $-parallel vector field which is known to solve
(11) on $M.$ By uniqueness, the unique solution $\xi ^{i}(x)$ is the
restriction of an $\varepsilon $-parallel vector field and therefore $\nabla
t=0$ implies that $t$ is $\varepsilon $-parallel if $t$ is a vector field.
In particular, we observe that we always have $\widetilde{\mathfrak{R}}=0$
on a parallelizable manifold $(M,\varepsilon ).$ Let $\mathfrak{X}(M)$
denote the Lie algebra of vector fields on $M$ and $\mathfrak{X}%
_{\varepsilon }(M)\subset \mathfrak{X}(M)$ denote the subspace of $%
\varepsilon $-invariant vector fields. We conclude that some $\xi \in 
\mathfrak{X}(M)$ belongs to $\mathfrak{X}_{\varepsilon }(M)$ if and only if
it solves (11) on $M.$ Now the integrability conditions of (10) for an
arbitrary tensor field $t$ is an expression in terms of $\widetilde{%
\mathfrak{R}}$ well known from tensor calculus. Therefore these conditions
are identically satisfied since $\widetilde{\mathfrak{R}}=0$ and we deduce
the desired implication for any tensor by a similar reasoning.

Clearly $\dim \mathfrak{X}_{\varepsilon }(M)=\dim M$. However $\mathfrak{X}%
_{\varepsilon }(M)$ need not be a Lie algebra, i.e., the bracket of two $%
\varepsilon $-parallel vector fields need not be $\varepsilon $-parallel. We
define

\begin{equation}
\widetilde{\nabla }_{r}\xi ^{i}\overset{def}{=}\frac{\partial \xi ^{i}}{%
\partial x^{r}}-\Gamma _{ar}^{i}\xi ^{a}
\end{equation}%
and extend $\widetilde{\nabla }$ to all tensor fields. Note that $r$ is now
the second index in $\Gamma _{\bullet r}^{\bullet }$ in (12). Assuming $%
\mathcal{R}=0,$ we check

\begin{equation}
\left( \frac{\partial \widetilde{\varepsilon }_{j}^{i}(x,y)}{\partial y^{k}}%
\right) _{y=x}=\Gamma _{jk}^{i}
\end{equation}%
So if we define $\widetilde{\Gamma }_{kj}^{i}$ by the LHS of (13) as in (6),
we get $\widetilde{\Gamma }_{kj}^{i}=\Gamma _{jk}^{i}.$ Recalling that $%
\widetilde{\varepsilon }$ is defined only if $\mathcal{R}=0,$ it is a
remarkable fact that $\widetilde{\nabla }$ is defined without the assumption 
$\mathcal{R}=0.$ If $\mathcal{R}=0,$ then $t$ is $\widetilde{\varepsilon }$%
-parallel (recall that this is a local condition) if and only if $\widetilde{%
\nabla }t=0.$

The integrability conditions of $\widetilde{\nabla }_{r}\xi ^{i}=0$ are
given by%
\begin{equation}
\mathfrak{R}_{rj,k}^{i}\overset{def}{=}\left[ \frac{\partial \Gamma _{kj}^{i}%
}{\partial x^{r}}+\Gamma _{kr}^{a}\Gamma _{aj}^{i}\right] _{[rj]}=0
\end{equation}

\begin{definition}
$\mathfrak{R}$ is the algebroid curvature of the parallelizable manifold $%
(M,\varepsilon ).$
\end{definition}

The following important proposition whose proof follows easily from
definitions (like all other facts in this section, except Proposition 4
below) clarifies the geometric meaning of $\mathfrak{R}$.

\begin{proposition}
$\mathfrak{X}_{\varepsilon }(M)$ is a Lie algebra if and only if $\mathfrak{R%
}=0.$ In this case, $\mathfrak{X}_{\widetilde{\varepsilon }}(M)$ is also a
Lie algebra ($\widetilde{\varepsilon }$ is defined since $\mathcal{R}=0$ by
Proposition 4) and the vector fields of $\mathfrak{X}_{\varepsilon }(M)$ and 
$\mathfrak{X}_{\widetilde{\varepsilon }}(M)$ commute.
\end{proposition}

Equation (12) is obtained from (2) by a "linearization" process whose
meaning will be clear shortly. In principle this process is the passage from
a groupoid to its algebroid. This formalism can be avoided in our simple
case of parallelizable manifolds but becomes indispensible for general
phg's. In the same way, $\mathfrak{R}$ is obtained from $\mathcal{R}$ by the
same linearization: substituting $y^{i}=x^{i}+t\xi ^{i}$ into $\mathcal{R}%
_{rj}^{m}(x,y)$ and differentiating with respect to $t$ at $t=0$ gives $%
\mathfrak{R}_{rj,a}^{m}\xi ^{a}.$ In particular, $\mathcal{R}=0$ implies $%
\mathfrak{R}=0.$

Now we have the following fundamental

\begin{proposition}
$\mathcal{R}=0$ $\Leftrightarrow \mathfrak{R}=0$
\end{proposition}

The implication $\Rightarrow $ states that "the Lie group has a Lie algebra"
whereas the nontrivial $\Leftarrow $ asserts that "the Lie algebra has a Lie
group" which is Lie's 3rd Fundamental Theorem.

Now suppose $\mathfrak{R}=0$ so that $(M,\varepsilon )$ is locally
homogeneous. The Lie algebra $\mathfrak{X}_{\widetilde{\varepsilon }}(M)$
integrates to the pseudogroup $\mathcal{S}$, i.e., $\mathfrak{X}_{\widetilde{%
\varepsilon }}(M)$ is the Lie algebra of the infinitesimal generators of the
action of $\mathcal{S}.$ Similarly, the Lie algebra $\mathfrak{X}%
_{\varepsilon }(M)$ integrates to the pseudogroup $\widetilde{\mathcal{S}}$,
in accordance with the familiar fact from Lie groups that the left (right)
invariant vector fields integrate to the right (left) actions. Recall,
however, that there is no canonical identification even if $\mathcal{S}$ is
globalizable.

Now we define the fundamental object

\begin{equation}
T_{jk}^{i}\overset{def}{=}\Gamma _{jk}^{i}-\Gamma _{kj}^{i}
\end{equation}

We have

\begin{equation}
\nabla _{r}\xi ^{i}=\widetilde{\nabla }_{r}\xi ^{i}+T_{ra}^{i}\xi ^{i}
\end{equation}%
and (16) easily generalizes to all tensor fields. The next proposition gives
the first hint that $T$ dominates the whole theory.

\begin{proposition}
\begin{equation}
\nabla _{r}T_{jk}^{i}=\mathfrak{R}_{jk,r}^{i}
\end{equation}
\end{proposition}

It follows that $\mathfrak{R}$ is determined by $T$ and $\mathfrak{R}=0$ if
and only if $T$ is $\varepsilon $-parallel!

To clarify the meaning of $T$ further, let $\xi ,\eta \in \mathfrak{X}(M).$
We define the torsion bracket $T(\xi ,\eta )\in \mathfrak{X}(M)$ by

\begin{equation}
T(\xi ,\eta )^{i}\overset{def}{=}T_{ab}^{i}\xi ^{a}\eta ^{b}
\end{equation}%
and the Jacobi $3$-form by

\begin{equation}
J(\xi ,\eta ,\sigma )\overset{def}{=}T(\xi ,T(\eta ,\sigma ))+T(\eta
,T(\sigma ,\xi ))+T(\sigma ,T(\xi ,\eta ))
\end{equation}

\begin{proposition}
(The First Bianchi Identity) Let $(M,\varepsilon )$ be parallelizable. Then
\end{proposition}

\begin{eqnarray}
&&\nabla _{\xi }T(\eta ,\sigma )+\nabla _{\eta }T(\sigma ,\xi )+\nabla _{\xi
}T(\eta ,\sigma )  \notag \\
&=&\mathfrak{R(}\eta ,\sigma )(\xi )+\mathfrak{R(}\sigma ,\xi )(\eta )+%
\mathfrak{R(}\xi ,\eta )(\sigma )  \notag \\
&=&J(\xi ,\eta ,\sigma )
\end{eqnarray}

In particular, if $\mathfrak{R}=0,$ then $J=0$. In this case, $T(\xi ,\eta
)=[\xi ,\eta ]$ for $\xi ,\eta \in \mathfrak{X}_{\varepsilon }(M),$ which
explains (15) \textit{to some extent}. For instance, certain odd degree
secondary characteristic classes are defined in [1] using $T$ with the
assumption $\mathfrak{R}=0.$

To finish this section, we recall that tensor calculus originated from
Riemannian geometry as an attempt to formalize Riemann's ideas. We hope to
have convinced the reader that tensor calculus (which we barely touched in
this section) could have originated also from Lie theory....and if this had
happened, then the concepts of torsion and curvature would have quite
different meanings today. We hope that the next section, where the above
formulas will be used in an essential way, will give further support to this
view.

\section{The homogeneous flow}

We recall the universal groupoid $\mathcal{U}_{1}$ and the subgroupoid $%
\varepsilon (\mathcal{U}_{0})=\varepsilon (M\times M)\subset \mathcal{U}%
_{1}. $ We fix some basepoint $\overline{\mathfrak{0}}\mathfrak{\in }M$ and
consider the principal bundle $\varepsilon (\overline{\mathfrak{0}}\times M)$
whose structure group is trivial as it is the $1$-arrow of the identity map
with source and target at $\overline{\mathfrak{0}}$. We fix some coordinates
around $\overline{\mathfrak{0}}$ once and for all.

We define a geometric object on $M$ whose components on $(U,x^{i})$ are $%
\varepsilon _{j}^{i}(\overline{\mathfrak{0}},x).$ Now (1) gives

\begin{equation}
\varepsilon _{a}^{i}(x,y)\varepsilon _{j}^{a}(\overline{\mathfrak{0}}%
,x)=\varepsilon _{j}^{i}(\overline{\mathfrak{0}},y)
\end{equation}

(21) asserts that $\varepsilon (M\times M)$ consists of those $1$-arrows in $%
\mathcal{U}_{1}$ which preserve the geometric object $\varepsilon _{j}^{i}(%
\overline{\mathfrak{0}},x).$ In view of (21), a change of coordinates $%
(U,x^{i})\rightarrow (V,y^{i})$ transforms $\varepsilon _{j}^{i}(\overline{%
\mathfrak{0}},x)$ by

\begin{equation}
\frac{\partial y^{i}}{\partial x^{a}}\varepsilon _{j}^{a}(\overline{%
\mathfrak{0}},x)=\varepsilon _{j}^{i}(\overline{\mathfrak{0}},y)
\end{equation}

(22) shows that $\varepsilon _{j}^{i}(\overline{\mathfrak{0}},y)$ transforms
only in the index $i$ but not in the index $j.$ We call $i$ the coordinate
index and $j$ the $\mathbb{R}^{n}$ index. We also define the dual object $%
\varepsilon _{j}^{i}(x,\overline{\mathfrak{0}})$ with $\varepsilon _{a}^{i}(%
\overline{\mathfrak{0}},x)\varepsilon _{j}^{a}(x,\overline{\mathfrak{0}}%
)=\varepsilon _{a}^{i}(x,\overline{\mathfrak{0}})\varepsilon _{j}^{a}(%
\overline{\mathfrak{0}}\mathfrak{,}x)=\delta _{j}^{i}$ where $i$ is the $%
\mathbb{R}^{n}$ index and $j$ is the coordinate index. If $\mathfrak{R}=0,$
it is an amusing fact to check that $\varepsilon _{j}^{i}(x,\overline{%
\mathfrak{0}})$ becomes the Maurer-Cartan form (see (51) in [1]).

Now differentiating (21) with respect to $y$ at $y=x$ and substituting from
(6) gives

\begin{equation}
\Gamma _{jk}^{i}(x)=\varepsilon _{k}^{a}(x,\overline{\mathfrak{0}})\frac{%
\varepsilon _{a}^{i}(\overline{\mathfrak{0}},x)}{\partial x^{j}}
\end{equation}

(23) shows that the RHS of (23) is independent of the choice of the base
point $\overline{\mathfrak{0}}.$ Now we identify $\overline{\mathfrak{0}}$
with the origin $\mathfrak{0}$ in $\mathbb{R}^{n}$ and (23) shows that we
can define $\Gamma _{jk}^{i}(x)$ consistently on the principal bundle $%
\varepsilon (\mathfrak{0}\times M).$ This identification will be useful in
Section 4. Note that the principal bundle $\varepsilon (\mathfrak{0}\times
M) $ determines the groupoid $\varepsilon (M\times M)$ since $\varepsilon
(p,q)=\varepsilon (\mathfrak{0},q)\circ \varepsilon (p,\mathfrak{0}).$

We rewrite (23) as

\begin{equation}
\nabla _{r}\varepsilon _{j}^{i}(\mathfrak{0},x)=\frac{\varepsilon _{j}^{i}(%
\mathfrak{0},x)}{\partial x^{r}}-\Gamma _{ra}^{i}(x)\varepsilon _{j}^{a}(%
\mathfrak{0},x)=0
\end{equation}

Similarly we have $\nabla _{r}\varepsilon _{j}^{i}(x,\mathfrak{0})=0$
keeping in mind that we always differentiate with respect to the coordinate
indices.

Now we define the geometric object $g$ by defining its components $g_{ij}(x)$
on $(U,x^{i})$ by

\begin{equation}
g_{ij}(x)\overset{def}{=}\varepsilon _{i}^{a}(x,\mathfrak{0})\varepsilon
_{j}^{a}(x,\mathfrak{0})=\varepsilon _{i}^{b}(x,\mathfrak{0})\left( \delta
_{ab}\varepsilon _{j}^{a}(x,\mathfrak{0})\right)
\end{equation}%
where we identify $\mathbb{R}^{n}$ with $\left( \mathbb{R}^{n}\right) ^{\ast
}$ by the canonical metric $\delta _{ij}.$ Clearly $g_{ij}$ is symmetric. It
is also positive definite since the matrix $\varepsilon _{j}^{i}(x,\mathfrak{%
0})$ is invertible. From (24) and (25) we deduce

\begin{equation}
\nabla _{r}g_{ij}=0
\end{equation}

\begin{definition}
$g$ is the canonical metric of the parallelizable manifold $(M,\varepsilon
). $
\end{definition}

Let $t_{km}^{ij}$ be a tensor field. We define

\begin{equation}
t_{km}^{i(j)}(x)\overset{def}{=}\varepsilon _{a}^{j}(x,\mathfrak{0}%
)t_{km}^{ia}(x)
\end{equation}

Now $t_{km}^{i(j)}(x)$ does not transform in the index $j$. We say that the
tensor $t_{km}^{i(j)}(x)$ is obtained from $t_{km}^{ij}(x)$ by moving the
index $j$ to $\mathfrak{0.}$ Similary we can move the index, say, $k$ to $%
\mathfrak{0}$ using $\varepsilon _{j}^{i}(\mathfrak{0},x)$ and this
operation extends to all tensors in an obvious way. With an abuse of
notation we will also move a covariant or contravariant $\mathbb{R}^{n}$
index $j$ to the coordinate index $(j)$ as in (33) below.

Now we define

\begin{eqnarray}
\mathfrak{H}_{j}^{i}(\varepsilon )\overset{def}{=}-\varepsilon _{j}^{a}(%
\mathfrak{0,}x)g^{bc}\nabla _{b}T_{ac}^{i} &=&-g^{bc}\nabla _{b}\varepsilon
_{j}^{a}(\mathfrak{0,}x)T_{ac}^{i}  \notag \\
&=&-g^{bc}\nabla _{b}T_{(j)c}^{i}  \notag \\
&=&-g^{bc}\mathfrak{R}_{(j)c,b}^{i}
\end{eqnarray}

Clearly $\mathfrak{R}=0$ implies $\mathfrak{H}=0.$ The converse will be
quite relevant in Section 5.

We now assume that the splitting $\varepsilon _{j}^{i}(\mathfrak{0},x,t)$
depends on time $t\geq 0$ smoothly and $\varepsilon (\mathfrak{0}%
,x,0)=\varepsilon _{0}.$ So for any small $t\geq 0$ and $x\in M,$ $%
\varepsilon (\mathfrak{0},x,t)$ assigns a $1$-arrow with source at $%
\mathfrak{0}$ and target at $x$ and this assignment is smooth in $x,t.$ We
observe that $\mathfrak{H}_{j}^{i}(\varepsilon )$ depends nonlinearly on the
second order derivatives of $\varepsilon .$ For simplicity of notation,
henceforth we omit the arguments of our functions, except those of $%
\varepsilon $ since the notation $\varepsilon $ does not distinguish between 
$\varepsilon (\mathfrak{0},x)$ and $\varepsilon (x,\mathfrak{0})$ which is
quite crucial below.

\begin{definition}
The homogeneous flow of a parallelizable manifold is the second order
nonlinear evolution equation

\begin{equation}
\frac{d\varepsilon _{j}^{i}(\mathfrak{0},x,t)}{dt}=\mathfrak{H}%
_{j}^{i}(\varepsilon )
\end{equation}%
with the initial condition $\varepsilon (\mathfrak{0},x,0)=\varepsilon _{0}.$
\end{definition}

Note that (29) stabilizes if $\mathfrak{H}=0$.

\begin{proposition}
If $M$ is compact, (29) admits a unique short time solution with any initial
condition.
\end{proposition}

\textbf{Proof : \ }We compute first the symbol of the linearization of $%
\mathfrak{H}$. So we set

\begin{equation}
\frac{d\varepsilon _{j}^{i}(\mathfrak{0},x,t)}{dt}=h_{j}^{i}
\end{equation}%
and compute the terms which depend on the second order derivatives of $%
h_{j}^{i}$ with respect $x$ in the expression

\begin{equation}
\frac{d\mathfrak{H}_{j}^{i}(\varepsilon (\mathfrak{0},x,t))}{dt}
\end{equation}

According to (24) $\nabla _{r}\varepsilon _{j}^{i}(\mathfrak{0},x,t)=0$ for
all $t$ where $\nabla $ is the operator defined by $\varepsilon (\mathfrak{0}%
,x,t)$ which we will write shortly as $\varepsilon (\mathfrak{0},x)$.
Therefore

\begin{eqnarray}
&&0=\frac{d}{dt}\nabla _{r}\varepsilon _{j}^{i}(\mathfrak{0},x)  \notag \\
&=&\frac{d}{dt}\left( \frac{\partial \varepsilon _{j}^{i}(\mathfrak{0},x)}{%
\partial x^{r}}-\Gamma _{ra}^{i}\varepsilon _{j}^{a}(\mathfrak{0},x)\right) 
\notag \\
&=&\frac{\partial }{\partial x^{r}}\left( \frac{d\varepsilon _{j}^{i}(%
\mathfrak{0},x)}{dt}\right) -\Gamma _{ra}^{i}\frac{d\varepsilon _{j}^{a}(%
\mathfrak{0},x)}{dt}-\frac{d\Gamma _{ra}^{i}}{dt}\varepsilon _{j}^{a}(%
\mathfrak{0},x)  \notag \\
&=&\frac{\partial h_{j}^{i}}{\partial x^{r}}-\Gamma _{ra}^{i}h_{j}^{a}-\frac{%
d\Gamma _{ra}^{i}}{dt}\varepsilon _{j}^{a}(\mathfrak{0},x)  \notag \\
&=&\nabla _{r}h_{j}^{i}-\frac{d\Gamma _{ra}^{i}}{dt}\varepsilon _{j}^{a}(%
\mathfrak{0},x)
\end{eqnarray}%
which gives

\begin{eqnarray}
\frac{d\Gamma _{rk}^{i}}{dt} &=&\varepsilon _{k}^{a}(x,\mathfrak{0})\nabla
_{r}h_{a}^{i}  \notag \\
&=&\nabla _{r}\left( \varepsilon _{k}^{a}(x,\mathfrak{0})h_{a}^{i}\right) 
\notag \\
&=&\nabla _{r}h_{(k)}^{i}
\end{eqnarray}

Note the simplicity of the variation (33) compared to the variation of the
Christoffel symbols of a metric in the Ricci flow.

From (33) we deduce

\begin{equation}
\frac{dT_{rk}^{i}}{dt}=\nabla _{r}h_{(k)}^{i}-\nabla _{k}h_{(r)}^{i}
\end{equation}

Now

\begin{eqnarray}
&&\frac{d\mathfrak{H}_{j}^{i}(\varepsilon )}{dt} \\
&=&-\frac{d\varepsilon _{j}^{a}(\mathfrak{0},x)}{dt}g^{bc}\nabla
_{b}T_{ac}^{i}-\varepsilon _{j}^{a}(\mathfrak{0},x)\frac{dg^{bc}}{dt}\nabla
_{b}T_{ac}^{i}-\varepsilon _{j}^{a}(\mathfrak{0},x)g^{bc}\frac{d}{dt}\nabla
_{b}T_{ac}^{i}  \notag
\end{eqnarray}

It is only the last term in (35) which contains second order derivatives of $%
h.$ Further

\begin{equation}
\frac{d}{dt}\nabla _{b}T_{ac}^{i}=\nabla _{b}\left( \frac{dT_{ac}^{i}}{dt}%
\right) +\text{ lower order terms}
\end{equation}

Substituting (36) into (35), the symbol is given by

\begin{eqnarray}
&&-\varepsilon _{j}^{a}(\mathfrak{0},x)g^{bc}\nabla _{b}\left( \frac{%
dT_{ac}^{i}}{dt}\right)  \notag \\
&=&-\varepsilon _{j}^{a}(\mathfrak{0},x)g^{bc}\nabla _{b}\left( \varepsilon
_{c}^{d}(x,\mathfrak{0})\nabla _{a}h_{d}^{i}-\varepsilon _{a}^{d}(x,%
\mathfrak{0})\nabla _{c}h_{d}^{i}\right)  \notag \\
&=&-\varepsilon _{j}^{a}(\mathfrak{0},x)\varepsilon _{c}^{d}(x,\mathfrak{0}%
)g^{bc}\nabla _{b}\nabla _{a}h_{d}^{i}+g^{bc}\nabla _{b}\nabla _{c}h_{j}^{i}
\end{eqnarray}%
and the second "elliptic" term in (37) shows that (29) is weakly parabolic.
\ 

Now we fix an arbitrary "connection" $\overline{\Gamma }_{jk}^{i}$ and
define the time-dependent vector field $W(x,t)$ by

\begin{equation}
W^{i}\overset{def}{=}g^{ab}\left( \Gamma _{ab}^{i}-\overline{\Gamma }%
_{ab}^{i}\right)
\end{equation}

The key fact in (38) is that $\Gamma -\overline{\Gamma }$ is a tensor and $%
\overline{\Gamma }$ does not depend on $t.$ We define the second order
nonlinear operator $\mathfrak{W}$ by the formula

\begin{equation}
\mathfrak{W}_{j}^{i}(\varepsilon )\overset{def}{=}\varepsilon _{j}^{a}(%
\mathfrak{0},x)\nabla _{a}W^{i}=\nabla _{(j)}W^{i}
\end{equation}

We compute

\begin{eqnarray}
\frac{d\mathfrak{W}_{j}^{i}}{dt} &=&\varepsilon _{j}^{c}(\mathfrak{0}%
,x)g^{ab}\nabla _{c}\frac{d\Gamma _{ab}^{i}}{dt}+\text{ }.....  \notag \\
&=&\varepsilon _{j}^{c}(\mathfrak{0},x)g^{ab}\nabla _{c}\left( \varepsilon
_{b}^{d}(x,\mathfrak{0})\nabla _{a}h_{d}^{i}\text{ }\right) +\text{ }..... 
\notag \\
&=&\varepsilon _{j}^{c}(\mathfrak{0},x)g^{ab}\varepsilon _{b}^{d}(x,%
\mathfrak{0})\nabla _{c}\nabla _{a}h_{d}^{i}+\text{ }.....  \notag \\
&=&\varepsilon _{j}^{c}(\mathfrak{0},x)g^{ab}\varepsilon _{b}^{d}(x,%
\mathfrak{0})\nabla _{a}\nabla _{c}h_{d}^{i}+\text{ }.....
\end{eqnarray}

From (37) and (40) we conclude that the evolution equation

\begin{equation}
\frac{d\varepsilon _{j}^{i}(\mathfrak{0},x)}{dt}=\mathfrak{H}%
_{j}^{i}(\varepsilon )+\mathfrak{W}_{j}^{i}(\varepsilon )
\end{equation}%
is strongly parabolic. By the well known existence theorem, we conclude that
(41) admits unique short time solutions. Now let $\varepsilon (\mathfrak{0}%
,x,t)=\varepsilon _{t}$ be the unique short time solution of (41) starting
from $\varepsilon _{0}.$ Let $\varphi _{t}$ be the unique short time
solution of the ODE

\begin{equation}
\frac{d\varphi _{t}}{dt}=W
\end{equation}%
so that $\varphi _{t}$ is a family of diffeomorphisms of $M$ with $\varphi
_{0}=id.$ Now it is easily shown (see [4]) that $\varphi _{t}^{\ast
}\varepsilon _{t}$ solves (29) and in fact the solutions of (29) are unique,
finishing the proof.

Note that the short time solution $\varepsilon _{j}^{i}(\mathfrak{0},x,t)$
is a $1$-arrow, i.e., an invertible matrix for \textit{sufficiently small }$%
t $ and for all $x\in M$ since this is so for the initial condition $%
\varepsilon _{j}^{i}(\mathfrak{0},x,0).$

\section{The gauge group}

In this section we will show that $\varepsilon (\mathfrak{0},x,t)$ depends
only on the value $\mathfrak{R}(x)$ for all $t$ in its domain by showing
that (29)  reduces to an ODE. For this purpose, we will first define the
first order universal group bundle $\mathcal{A\rightarrow }$ $M.$

Let $\mathcal{U}_{1}^{p,p}$ denote the set of all $1$-arrows with source and
target at $p\in M.$ A choice of coordinates around $p$ identifies $\mathcal{U%
}_{1}^{p,p}$ with the Lie group $GL(n,\mathbb{R}).$ We define $\mathcal{A}%
\overset{def}{=}\cup _{p\in M}\mathcal{U}_{1}^{p,p}$. A local section of $%
\mathcal{A}$ is of the form $\mathfrak{a}_{j}^{i}(x)$ in coordinates and is
an invertible linear map on the tangent space at $x.$ Let $\Gamma \mathcal{A}
$ denote the space of smooth sections of $\mathcal{A}$. Now $\Gamma \mathcal{%
A}$ is a group with the fiberwise composition of jets. We call elements of $%
\Gamma \mathcal{A}$ (first order) gauge transformations.

Now let $M$ be a parallelizable manifold and $\mathcal{E}$ denote the set of
all splittings on $M.$ We recall that for any $x\in M$ the splitting $%
\varepsilon $ assigns a $1$-arrow from $\mathfrak{0\in }\mathbb{R}^{n}$ to $%
x $ and we also write $\varepsilon (\mathfrak{0}\times M)$ for $\varepsilon
. $ Let $\mathfrak{g}\in \Gamma \mathcal{A}$ and $\varepsilon (\mathfrak{0}%
\times M)\in \mathcal{E}$. We define $\mathfrak{g}\varepsilon \in \mathcal{E}
$ by

\begin{equation}
\left( \mathfrak{a}\varepsilon \right) (\mathfrak{0},p)\overset{def}{=}%
\mathfrak{a}(p)\circ \varepsilon (\mathfrak{0},p)\text{ \ \ \ }p\in M
\end{equation}

The action (43) is easily seen to be simply transitive. Hence a choice of $%
\varepsilon $ gives a 1-1 correspondence between $\Gamma \mathcal{A}$ and $%
\mathcal{E}.$ A straightforward computation using (43), (6) and (14) gives
the well known formula

\begin{equation}
\mathfrak{R}\left( \mathfrak{a}\varepsilon \right) _{jk,m}^{i}=\mathfrak{a}%
_{a}^{i}\mathfrak{R}\left( \varepsilon \right) _{jk,b}^{a}\mathfrak{b}%
_{m}^{b}\text{ \ \ \ \ \ \ }\mathfrak{a}_{a}^{i}\mathfrak{b}_{j}^{a}=\delta
_{j}^{i}
\end{equation}

Now $\Gamma \mathcal{A}$ is an infinite dimensional Frechet Lie group as
follows (see [7] for the technical detais of this theory). Let $\xi $ be a
section of the vector bundle $Hom(T,T)\rightarrow M,$ so that $\xi (p)$ is a
linear map at $T_{p}.$ We note that the Lie algebra of $\mathcal{U}%
_{1}^{p,p} $ is canonically isomorphic to the fiber $Hom(T_{p},T_{p})$
endowed with the usual bracket of matrices. Therefore the space of sections $%
\Gamma Hom(T,T)$ is a Lie algebra. Note that the bracket of $\Gamma Hom(T,T)$
is defined pointwise and does not involve differentiation. Now the Lie
algebra $\Gamma Hom(T,T)$ may be viewed (and this view can be made rigorous)
as the Lie algebra of $\Gamma \mathcal{A}$ as follows. The value $\xi (p)$
determines a $1$-parameter subgroup $\gamma (p,t)$ in $\mathcal{U}_{1}^{p,p}$
with $\left( \frac{d\gamma (p,t)}{dt}\right) _{t=0}=\xi (p),$ $\gamma
(p,0)=p $ for all $p\in M$ and $\gamma (p,t)$ is defined for all $t\geq 0.$
So $\xi \in \Gamma Hom(T,T)$ determines a "$1$-parameter subgroup" of $%
\Gamma \mathcal{A} $ which consists of the $1$-parameter subgroups at all
points. The local bijectivity of the pointwise exponential maps gives the
local bijectivity of the exponential map $\exp :\Gamma Hom(T,T)\rightarrow
\Gamma \mathcal{A}$ if $M$ is compact.

Now let $\varepsilon (\mathfrak{0},x,t)$ be the unique short time solution
of (29) starting from $\varepsilon (\mathfrak{0},x,0)=\varepsilon _{0}.$ In
view of the simple transitivity of (43), for any small $t\geq 0$ there
exists a unique $\mathfrak{a}\in \Gamma \mathcal{A}$ satisfying $\varepsilon
(\mathfrak{0},x,t)=\mathfrak{a}(x,t)\circ \varepsilon (\mathfrak{0},x,0)$
where $\mathfrak{a}(x,0)=id.$ In coordinates

\begin{equation}
\varepsilon _{j}^{i}(\mathfrak{0},x,t)=\mathfrak{a}_{a}^{i}(x,t)\varepsilon
_{j}^{a}(\mathfrak{0},x,0)
\end{equation}

Using (25) and (45) we find

\begin{equation}
g^{ij}(x,t)=g^{ab}(x,0)\mathfrak{a}_{a}^{i}(x,t)\mathfrak{a}_{b}^{j}(x,t)
\end{equation}

We now substitute (45) into (29) and using (44) and (46) we compute

\begin{eqnarray}
\frac{d\varepsilon _{j}^{i}(\mathfrak{0},x,t)}{dt} &=&\frac{d\mathfrak{a}%
_{a}^{i}(x,t)}{dt}\varepsilon _{j}^{a}(\mathfrak{0},x,0)  \notag \\
&=&\mathfrak{H}_{j}^{i}(\mathfrak{a}_{t}\varepsilon _{0})  \notag \\
&=&\mathfrak{R}_{(j)a,b}^{i}\mathfrak{(\mathfrak{a}_{t}\varepsilon _{0})}%
g^{ab}(x,t)  \notag \\
&=&\varepsilon _{j}^{e}(\mathfrak{0},x,t)\mathfrak{R}_{ea,b}^{i}(\mathfrak{%
\mathfrak{a}_{t}\varepsilon _{0}})g^{cd}(x,0)\mathfrak{a}_{c}^{a}(x,t)%
\mathfrak{a}_{d}^{b}(x,t)  \notag \\
&=&\mathfrak{a}_{f}^{e}(x,t)\varepsilon _{j}^{f}(\mathfrak{0},x,0)\mathfrak{a%
}_{g}^{i}(x,t)\mathfrak{R}_{ea,h}^{g}(\mathfrak{\varepsilon _{0}})\mathfrak{b%
}_{b}^{h}(x,t)g^{cd}(x,0)\mathfrak{a}_{c}^{a}(x,t)\mathfrak{a}_{d}^{b}(x,t) 
\notag \\
&=&\mathfrak{a}_{f}^{e}(x,t)\varepsilon _{j}^{f}(\mathfrak{0},x,0)\mathfrak{a%
}_{g}^{i}(x,t)\mathfrak{R}_{ea,d}^{g}(\mathfrak{\varepsilon _{0}})g^{cd}(x,0)%
\mathfrak{a}_{c}^{a}(x,t)
\end{eqnarray}

Therefore

\begin{equation}
\frac{d\mathfrak{a}_{j}^{i}(x,t)}{dt}=\mathfrak{a}_{j}^{e}(x,t)\mathfrak{a}%
_{g}^{i}(x,t)\mathfrak{R}_{ea,d}^{g}(\mathfrak{\varepsilon _{0}})g^{cd}(x,0)%
\mathfrak{a}_{c}^{a}(x,t)
\end{equation}

The equation (48) is an ODE for fixed $x$ which has a unique analytic
solution $\mathfrak{a}_{j}^{i}(x,t)$ for $-\epsilon \leq t\leq \epsilon $
with the inital condition $\mathfrak{a}_{j}^{i}(x,0)=\delta _{j}^{i}.$

Thus we proved

\begin{proposition}
$\varepsilon (\mathfrak{0},x,t)$ is the unique short time solution of (29)
with the initial condition $\varepsilon _{0}$ if and only if $\mathfrak{a}%
(x,t)$ defined by (45) is the unique solution of the ODE (48).
\end{proposition}

An inspection of (48) indicates that the power series determined by (48)
centered at $t=0$ will not converge for all $t\in \mathbb{R}$ in general.
This fact can be most easily seen by checking the radius of convergence of
the power series at $t=0$ of the 1-dimensional analog $\frac{d\mathfrak{a(t)}%
}{dt}=\mathfrak{a}^{3}(t)R$ of (48). Thus HF will develop finite time
singularities like the Ricci flow.

Note that (48) gives

\begin{eqnarray}
\frac{d\mathfrak{a}_{j}^{i}(x,0)}{dt} &=&\mathfrak{R}_{ja,b}^{i}(\mathfrak{%
\varepsilon _{0}})g^{ab}(x,0)  \notag \\
&=&\mathfrak{H}_{(j)}^{i}(\varepsilon _{0})
\end{eqnarray}%
and the 1-parameter subgroup of $\Gamma \mathcal{A}$ determined by the
initial condition (49) is defined by

\begin{equation}
\frac{d\mathfrak{a}_{j}^{i}(x,t)}{dt}=\mathfrak{a}_{a}^{i}(x,t)\mathfrak{H}%
_{(j)}^{a}(\varepsilon _{0})
\end{equation}%
whose power series converges for all $t\in \mathbb{R}.$

\section{Poincare Conjecture}

Assume that $M$ is a parallelizable manifold. We recall here that an
orientable 3-manifold is parallelizable. Consider the following two
assertions.

\textbf{A1}. Suppose $\varepsilon _{t}$ converges to some parallelism $%
\varepsilon _{\infty }$ as $t\rightarrow \infty ,$ i.e., no singularities
develop including $t=\infty .$ Then $\mathfrak{H(}\varepsilon _{\infty })=0$.

\textbf{A2. }If $\dim M=3,$ then $\mathfrak{H}=0$ $\Leftrightarrow \mathfrak{%
R}=0.$ Further, the hypothesis of \textbf{A1} holds for some $\varepsilon
_{0}$ if $M$ is compact and simply connected.

\begin{proposition}
\textbf{A1 }and \textbf{A2} imply PC.
\end{proposition}

\textbf{Proof : }It suffices to show that the compact and simply connected
local Lie group $(M,\varepsilon _{\infty })$ of dimension three is
diffeomorphic to $S^{3}.$ Since $M$ is compact, $\varepsilon _{\infty }$ is
easily seen to be complete (see Lemma 7.3 in [1]). Since $M$ is also simply
connected, the pseudogroup $\mathcal{S}$ in Section 2 globalizes to a Lie
group. However $S^{3}$ is the only compact and simply connected Lie group in
dimension three up to diffeomorphism.

\section{ Appendix. The relation of HF to the Ricci flow}

Suppose we evolve the initial canonical metric given by Definition 7
according to the Ricci flow

\begin{equation}
\frac{dg}{dt}=-2Ric(g)
\end{equation}

The natural question is how (29) and (51) are related. So let $g_{t}$ be the
unique short time solution of (51) with $g_{0}$ being the canonical metric
of $\varepsilon _{0}$ so that $\varepsilon _{0}(\mathfrak{0},M)$ is a
trivialization of the $O(n)$-principal bundle $P(g_{0})$ determined by $%
g_{0}.$ It is easy to see that the principal bundles $P(g_{t})$ admit
trivializations for small $t.$ However, there is no canonical way of
choosing these trivializations. Therefore, (51) does not imply (29). The
main idea of this note is to reverse this reasoning and try to derive (51)
from (29) as follows. Now (26) gives

\begin{equation}
\frac{\partial g_{ij}}{\partial x^{r}}+\Gamma _{ri}^{a}g_{aj}+\Gamma
_{rj}^{a}g_{ia}=0
\end{equation}

Let $\Sigma _{jk}^{i}$ be the Christoffel symbols of $g_{ij}$ so that (52)
holds also for $\Sigma _{jk}^{i}.$ The Gauss trick of shifting the indices
in (52) gives the formula

\begin{equation}
\Sigma _{jk}^{i}=-\frac{1}{2}\left( \Gamma
_{jk}^{i}+T_{jb}^{a}g_{ka}g^{ib}\right) _{(jk)}
\end{equation}%
where $(jk)$ denotes the symmetrization of $j,k.$ It is natural to
substitute (53) into (51) and try to derive (51) from (29). Equivalently, we
may differentiate (25), substitute from (29) and try to express the
resulting expression in terms of $g$ by eliminating $\varepsilon $. Unable
to carry out this derivation, we came up with the expression (28). We now
believe that (29) and (51) are independent.

\bigskip

\bigskip \bigskip

\bigskip

\bigskip

\bigskip

\textbf{References}

\bigskip

\bigskip \lbrack 1] \ E.Abado\u{g}lu, E.Orta\c{c}gil: Intrinsic
characteristic classes of a local Lie group, Portugal. Math. (N.S.), Vol.67,
Fasc.4, 2010, 453-483

[2] A.D.Blaom: Lie algebroids and Cartan's method of equivalence, Trans.
Amer. Math. Soc., 364, (2012), 3071-3135

[3] D.DeTurck: Deforming metrics in the direction of their Ricci tensors, J.
Diff. Geometry, 18, 157-162, 1983

[4] R.Hamilton: Three manifolds with positive Ricci curvature, J. Diff.
Geometry, 17, 255-306, 1982

[5] P.J.Olver: Nonassociative local Lie groups, J. Lie Theory, 6, (1996),
23-51

[6] P.J.Olver: private communication

[7] H.Omori: Infinite dimensional Lie groups, Translations of Mathematics
Monographs, Vol. 128, AMS, 1996

[8] E.Orta\c{c}gil: Characteristic classes as obstructions to local
homogeneity, arXiv: 1401,1116, extended version of the talk given at the
Fields Institute, Toronto, Focused Workshop on Exterior Differential Systems
and Lie Theory, Dec. 9-13, 2013, organized by R.Fernandes, N.Kamran,
P.J.Olver

\bigskip

\bigskip

Erc\"{u}ment Orta\c{c}gil

\bigskip

Bogazici University (Emeritus), Mathematics Department, Bebek, 34342,
Istanbul, Turkey

now at: Bodrum, Mu\u{g}la, Turkey \ \ \ \ \ \ \ \ \ \ \ 

\bigskip

ortacgil@boun.edu.tr

ortacgile@gmail.com

\end{document}